\documentclass[draft,12pt,reqno]{amsart}
\usepackage{amsmath,amsthm,amssymb,amscd}

\begin{document}

\newcommand{\DATE}{July 2004}
\newcommand{\TITLE}{Generalized Greatest Common Divisors,
Divisibility Sequences, and Vojta's Conjecture for Blowups}
\newcommand{\TITLERUNNING}{Generalized GCD's and Vojta's Conjecture for Blowups}
\newcommand{\AUTHOR}{Joseph H. Silverman}

\theoremstyle{plain}
\newtheorem{theorem}{Theorem}
\newtheorem{conjecture}[theorem]{Conjecture}
\newtheorem{proposition}[theorem]{Proposition}
\newtheorem{lemma}[theorem]{Lemma}
\newtheorem{corollary}[theorem]{Corollary}

\theoremstyle{definition}
\newtheorem{definition}{Definition}

\theoremstyle{remark}
\newtheorem{remark}{Remark}
\newtheorem{example}{Example}
\newtheorem{question}{Question}
\newtheorem*{acknowledgement}{Acknowledgements}

\def\BigStrut{\vphantom{$(^{(^(}_{(}$}} 

\newenvironment{notation}[0]{%
  \begin{list}%
    {}%
    {\setlength{\itemindent}{0pt}
     \setlength{\labelwidth}{4\parindent}
     \setlength{\labelsep}{\parindent}
     \setlength{\leftmargin}{5\parindent}
     \setlength{\itemsep}{0pt}
     }%
   }%
  {\end{list}}

\newenvironment{parts}[0]{%
  \begin{list}{}%
    {\setlength{\itemindent}{0pt}
     \setlength{\labelwidth}{1.5\parindent}
     \setlength{\labelsep}{.5\parindent}
     \setlength{\leftmargin}{2\parindent}
     \setlength{\itemsep}{0pt}
     }%
   }%
  {\end{list}}
\newcommand{\Part}[1]{\item[\upshape#1]}

\renewcommand{\a}{\alpha}
\renewcommand{\b}{\beta}
\newcommand{\g}{\gamma}
\renewcommand{\d}{\delta}
\newcommand{\e}{\epsilon}
\newcommand{\f}{\phi}
\renewcommand{\l}{\lambda}
\renewcommand{\k}{\kappa}
\newcommand{\lhat}{\hat\lambda}
\newcommand{\m}{\mu}
\renewcommand{\o}{\omega}
\newcommand{\p}{\pi}
\renewcommand{\r}{\rho}
\newcommand{\rbar}{{\bar\rho}}
\newcommand{\s}{\sigma}
\newcommand{\sbar}{{\bar\sigma}}
\renewcommand{\t}{\tau}
\newcommand{\z}{\zeta}

\newcommand{\D}{\Delta}
\newcommand{\F}{\Phi}
\newcommand{\G}{\Gamma}

\newcommand{\ga}{{\mathfrak{a}}}
\newcommand{\gb}{{\mathfrak{b}}}
\newcommand{\gc}{{\mathfrak{c}}}
\newcommand{\gd}{{\mathfrak{d}}}
\newcommand{\gm}{{\mathfrak{m}}}
\newcommand{\gn}{{\mathfrak{n}}}
\newcommand{\gp}{{\mathfrak{p}}}
\newcommand{\gq}{{\mathfrak{q}}}
\newcommand{\gP}{{\mathfrak{P}}}
\newcommand{\gQ}{{\mathfrak{Q}}}

\def\Acal{{\mathcal A}}
\def\Bcal{{\mathcal B}}
\def\Ccal{{\mathcal C}}
\def\Dcal{{\mathcal D}}
\def\Ecal{{\mathcal E}}
\def\Fcal{{\mathcal F}}
\def\Gcal{{\mathcal G}}
\def\Hcal{{\mathcal H}}
\def\Ical{{\mathcal I}}
\def\Jcal{{\mathcal J}}
\def\Kcal{{\mathcal K}}
\def\Lcal{{\mathcal L}}
\def\Mcal{{\mathcal M}}
\def\Ncal{{\mathcal N}}
\def\Ocal{{\mathcal O}}
\def\Pcal{{\mathcal P}}
\def\Qcal{{\mathcal Q}}
\def\Rcal{{\mathcal R}}
\def\Scal{{\mathcal S}}
\def\Tcal{{\mathcal T}}
\def\Ucal{{\mathcal U}}
\def\Vcal{{\mathcal V}}
\def\Wcal{{\mathcal W}}
\def\Xcal{{\mathcal X}}
\def\Ycal{{\mathcal Y}}
\def\Zcal{{\mathcal Z}}

\renewcommand{\AA}{\mathbb{A}}
\newcommand{\BB}{\mathbb{B}}
\newcommand{\CC}{\mathbb{C}}
\newcommand{\FF}{\mathbb{F}}
\newcommand{\GG}{\mathbb{G}}
\newcommand{\NN}{\mathbb{N}}
\newcommand{\PP}{\mathbb{P}}
\newcommand{\QQ}{\mathbb{Q}}
\newcommand{\RR}{\mathbb{R}}
\newcommand{\ZZ}{\mathbb{Z}}

\def \bfa{{\mathbf a}}
\def \bfb{{\mathbf b}}
\def \bfc{{\mathbf c}}
\def \bfe{{\mathbf e}}
\def \bff{{\mathbf f}}
\def \bfF{{\mathbf F}}
\def \bfg{{\mathbf g}}
\def \bfp{{\mathbf p}}
\def \bfr{{\mathbf r}}
\def \bfs{{\mathbf s}}
\def \bft{{\mathbf t}}
\def \bfu{{\mathbf u}}
\def \bfv{{\mathbf v}}
\def \bfw{{\mathbf w}}
\def \bfx{{\mathbf x}}
\def \bfy{{\mathbf y}}
\def \bfz{{\mathbf z}}


\newcommand{\Kbar}{{\bar K}}
\newcommand{\Vbar}{{\bar V}}
\newcommand{\Wbar}{{\bar W}}

\newcommand{\Atilde}{{\tilde A}}
\newcommand{\Ptilde}{{\tilde P}}
\newcommand{\Vtilde}{{\tilde V}}
\newcommand{\Wtilde}{{\tilde W}}
\newcommand{\Xtilde}{{\tilde X}}
\newcommand{\Ytilde}{{\tilde Y}}
\newcommand{\Ztilde}{{\tilde Z}}

\newcommand{\Ghat}{{\hat \Gamma}}
\newcommand{\Xhat}{{\hat X}}


\newcommand{\Aut}{\operatorname{Aut}}
\newcommand{\bdy}{\partial}
\newcommand{\Disc}{\operatorname{Disc}}
\newcommand{\Div}{\operatorname{Div}}
\newcommand{\End}{\operatorname{End}}
\newcommand{\Frob}{\operatorname{Frob}}
\newcommand{\Gal}{\operatorname{Gal}}
\newcommand{\Gbar}{{\bar G}}
\newcommand{\GCD}{\operatorname{GCD}}
\newcommand{\Gtilde}{{\tilde G}}
\newcommand{\hhat}{{\hat h}}
\newcommand{\Ideal}{\operatorname{Ideal}}
\newcommand{\Image}{\operatorname{Image}}
\newcommand{\longhookrightarrow}{\lhook\joinrel\relbar\joinrel\rightarrow}
\newcommand{\LS}[2]{\genfrac(){}{}{#1}{#2}}  
\newcommand{\MOD}[1]{~(\textup{mod}~#1)}
\newcommand{\Norm}{\operatorname{N}}
\newcommand{\notdivide}{\nmid}
\newcommand{\ord}{\operatorname{ord}}
\newcommand{\Pic}{\operatorname{Pic}}
\newcommand{\Proj}{\operatorname{Proj}}
\newcommand{\Qbar}{{\bar{\QQ}}}
\newcommand{\rank}{\operatorname{rank}}
\newcommand{\res}{\operatornamewithlimits{res}}
\newcommand{\Resultant}{\operatorname{Resultant}}
\renewcommand{\setminus}{\smallsetminus}
\newcommand{\Spec}{\operatorname{Spec}}
\newcommand{\Support}{\operatorname{Support}}
\newcommand{\tors}{{\textup{tors}}}
\newcommand{\<}{\langle}
\renewcommand{\>}{\rangle}


\title[\TITLERUNNING]{\TITLE}
\date{\DATE}
\author{\AUTHOR}
\address{Mathematics Department, Box 1917, Brown University,
Providence, RI 02912 USA}
\email{jhs@math.brown.edu}
\thanks{The author's research supported by NSA grant H98230-04-1-0064}

\begin{abstract}
We apply Vojta's conjecture to blowups and deduce a
number of deep statements regarding (generalized) greatest common
divisors on varieties, in particular on projective space and on 
abelian varieties. Special cases of these statements generalize
earlier results and conjectures. We also discuss the relationship
between generalized greatest common divisors and
the divisibility sequences attached to algebraic groups, and we
apply Vojta's conjecture to obtain a strong bound on
the divisibility sequences attached to abelian varieties of dimension
at least two.
\end{abstract}

\maketitle

\section*{Introduction}
Bugeaud, Corvaja, and Zannier~\cite{BCZ} recently proved that if~$a$
and~$b$ are multiplicatively independent integers, then for every
$\e>0$ there is a constant $N=N(a,b,\e)$ so that
\begin{equation}
  \label{equation:BCZ}
  \gcd(a^n-1,b^n-1) \le 2^{\e n}
  \qquad\text{for all $n\ge N$.}
\end{equation}
The proof of this beautiful, but innocuous looking, inequality
requires an ingenious application of Schmidt's Subspace
Theorem~\cite{Schmidt}.
Corvaja and Zannier~\cite[Proposition~4]{CZ2}
generalize~\eqref{equation:BCZ} by replacing~$a^n$ and~$b^n$ with
arbitrary elements from a fixed finitely generated subgroup
of~$\Qbar^*$.  
For ease of exposition, we state their result over~$\QQ$.

\begin{theorem}[Corvaja-Zannier \cite{CZ2}]
\label{theorem:CZunitminusone}
Let $S$ be a finite set of rational primes and let $\e>0$.
There is a finite set $Z=Z(S,\e)\subset\ZZ^2$
so that all $\a,\b\in\ZZ_S^*\cap\ZZ$
satisfy one of
the following three conditions:
\begin{parts}
\item[\upshape(1)]
$(\a,\b)\in Z$.
\item[\upshape(2)]
$\a^m=\b^n$ for some $(m,n)$ satisfying $1\le\max\{m,n\}\le\e^{-1}$.
\item[\upshape(3)]
$\gcd(\a-1,\b-1) \le  \max\bigl(|\a|,|\b|\bigr)^\e$.
\end{parts}
\end{theorem}

In other words, if~$\a,\b\in\ZZ$ are $S$-units, then
\[
  \gcd(\a-1,\b-1) \le \max\bigl(|\a|,|\b|\bigr)^\e
\]
except for some obvious families of exceptions together with a
finite number of additional exceptions.  Analogous statements for
elliptic curves and/or over function fields have been studied by a
number of authors \cite{AR,McKinnon,JHSFF,JHSECFF}.
\par
The purpose of this note is to explain how Vojta's
Conjecture~\cite[Conjecture~3.4.3]{Vojta} applied to varieties blown
up along smooth subvarieties leads to a very general statement about
greatest common divisors that encompasses many known results and
previous conjectures. Thus although we do not prove 
unconditional results in this paper, we hope that the application of Vojta's
conjecture will help to put the problem of gcd bounds
into a general context, while at the same time suggesting precise
statements whose proofs may be possible using current techniques from
Diophantine approximation and arithmetic geometry. (See also
McKinnon's paper~\cite{McKinnon} for a discussion of Vojta's conjecture applied
to certain blowups.)
\par
We begin in the
Section~\ref{section:examplesoverQ}
by describing three special cases of our main theorem. These
serve to motivate our general result and to justify the notation that
is needed later. 
We next in Section~\ref{section:generalizedgcddef}
set notation and explain how a generalized concept of
greatest common divisor is naturally formulated in terms of the
height of points on blowup varieties with respect to the exceptional
divisor of the blowup.
Section~\ref{section:vojtaconjecture}
states Vojta's conjecture, followed
in Section~\ref{section:applyvojtatoblowups}
by our main result (Theorem~\ref{theorem:vojta=>hgcdbound})
in which we apply Vojta's conjecture to a blowup variety, making use
of the well-known relation between the canonical bundle on a variety
and on its blowup. In Section~\ref{section:proofsof3theorems}
we apply our main theorem to prove the three special cases from
Section~\ref{section:examplesoverQ}, including some additional arguments to
pin down the exceptional sets more precisely. 
Section~\ref{section:divseqsandalggps}
takes up the question of divisibility sequences, which are
sequences~$(a_n)_{n\ge1}$ satisfying \text{$m|n\Rightarrow a_m|a_n$}.
We are especially interested in divisibility sequences associated to
algebraic groups, or more precisely, to group schemes over~$\ZZ$.
We show that these geometric divisibility sequences are closely related
to generalized greatest common divisors and apply Vojta's conjecture
to the divisibility sequences attached to abelian varieties of dimension
at least~2. Finally, in Section~\ref{section:finalremarks},
we make a few final remarks and pose some questions.

\begin{acknowledgement}
I would like to thank
Y. Bugeaud,
P. Corvaja,
D. McKinnon,
Z. Rudnick,
G. Walsh,
and
U. Zannier,
for their helpful email correspondence during the preparation of 
this article, and D. Abramovich,
P. Corvaja, and U. Zannier
for correcting errors in the initial draft.
\end{acknowledgement}

\section{Three special cases over $\QQ$}
\label{section:examplesoverQ}

In this section we describe three special cases of our main
theorem. These generalize earlier results and conjectures appearing in
the literature.  In order to avoid excessive notation, we restrict
ourselves to working over~$\QQ$. All results are conditional on the
validity of Vojta's conjeture.  We refer the reader to
Section~\ref{section:vojtaconjecture}
(Conjecture~\ref{conjecture:vojta}) or to Vojta's original
monograph~\cite[Conjecture~3.4.3]{Vojta} for the statement of Vojta's
conjecture.  In order to state our first result, we need one piece of
notation.

\begin{definition}
Let $S$ be a finite set of rational primes. For any nonzero
integer~$x\in\ZZ$, we write
$|x|_S'$ for the largest divisor of~$x$ that is not divisible by any of
the primes in~$S$, i.e. 
\[
  |x|_S'=|x|\prod_{p\in S}|x|_p.
\]
Informally, we call~$|x|_S'$ the ``prime-to-$S$'' part of~$x$. In
particular,~$x$ is an $S$-unit if and only if $|x|_S'=1$.
\end{definition}

Our first result deals with~$\PP^n$ blown up along a smooth subvariety.

\begin{theorem}
\label{theorem:PnoverQ}
Fix a finite set of rational primes~$S$.
Let $f_1,f_2,\ldots,f_t\in\ZZ[X_0,\ldots,X_n]$ be homogeneous polynomials
so that the set of zeros
\[
  V = \{f_1=f_2=\cdots=f_t=0\} \subset \PP^n
\]
is a smooth variety, and assume further that~$V$ does not intersect the union
of the coordinate hyperplanes~$\bigcup_{i=0}^n\{X_i=0\}$.  Let
$r=n-\dim(V)\ge2$ be the codimension of~$V$ in~$\PP^n$.
\par
Assume that Vojta's conjecture is true {\upshape(}for~$\PP^n$ blown up
along~$V${\upshape)}.  Fix $\e>0$. Then there is a homogeneous
polynomial $0\ne g\in\ZZ[X_0,\ldots,X_n]$, depending
on~$f_1,\ldots,f_t$ and~$\e$, and a constant~$\d$, depending
on~$f_1,\ldots,f_t$, so that every
$\bfx=(x_0,x_1,\ldots,x_n)\in\ZZ^{n+1}$ with $\gcd(x_0,\ldots,x_n)=1$
satisfies either
\begin{parts}
\item[\upshape(1)]
$g(\bfx)=0$, or 
\item[\upshape(2)]
$\gcd\bigl(f_1(\bfx),\ldots,f_t(\bfx)\bigr)$ 
\hfil\break
\phantom{$\gcd\bigl(f_1(\bfx),\ldots,$}
${}\le \max\bigl\{|x_0|,\ldots,|x_n|\bigr\}^\e
\cdot\left(|x_0x_1\cdots x_n|_S'\right)^{1/(r-1+\d\e)}.$
\end{parts}
\end{theorem}

\begin{example}
We apply Theorem~\ref{theorem:PnoverQ} to~$\PP^2$ with $f_1=X_1-X_0$
and $f_2=X_2-X_0$. Then~$V$ is a single point and $r=2$, so the
theorem says that off of a one dimensional exceptional set we have
\begin{equation} 
  \label{equation:gcdx0x1x2}
  \gcd(x_1-x_0,x_2-x_0) \le \max\bigl\{|x_0|,|x_1|,|x_2|\bigr\}^\e
   \cdot \left(|x_0x_1x_2|_S'\right)^{1/(1+\d\e)}.
\end{equation} 
In particular, suppose that we take $x_0=1$ and restrict~$x_1$ and~$x_2$
to be $S$-units, as in Theorem~\ref{theorem:CZunitminusone}.
Then $|x_0x_1x_2|_S'=1$, so~\eqref{equation:gcdx0x1x2} becomes
\[
  \gcd(x_1-1,x_2-1) \le \max\bigl\{|x_1|,|x_2|\bigr\}^\e
\]
and we recover the inequality of Theorem~\ref{theorem:CZunitminusone}.
(Theorem~\ref{theorem:CZunitminusone} also includes a description of
the exceptional set, but once one knows that the exceptional set is a
union of curves, it is not hard to recover this description.) Thus
Vojta's conjecture implies a natural generalziation of
Theorem~\ref{theorem:CZunitminusone} in which we remove the
restriction that~$\a$ and~$\b$ be $S$-units and replace condition~(3)
of Theorem~\ref{theorem:CZunitminusone} with the inequality
\begin{equation}
  \label{equation:gcda-1b-1}
  \gcd(\a-1,\b-1) \le \max\bigl\{|\a|,|\b|\bigr\}^\e
   \cdot \left(|\a\b|_S'\right)^{1/(1+\d\e)}.
\end{equation}
It would be quite interesting to give an unconditional proof of this
generalization. We also remark that a closer analysis of this special
case of Theorem~\ref{theorem:PnoverQ} shows
that~\eqref{equation:gcda-1b-1} should be valid for any \text{$\d<1$ }.
\end{example}

Our second example deals with elliptic  curves and has applications to
the theory of elliptic divisibility sequences.

\begin{theorem}
\label{theorem:ECoverQ}
Let $E/\QQ$ be an elliptic curve given by a Weierstrass equation, and
for any nonzero point $P=(x_P,y_P)\in E(\QQ)$, write $x_P=A_P/D_P^2$
as a fraction in lowest terms with $D_P>0$. Also let
$H(P)=H(x_P)=\max\{|A_P|,|D_P^2|\}$ be the usual Weil height on~$E$.
\par
Assume that Vojta's conjecture is true {\upshape(}for~$E^2$ blown up
at~$(O,O)${\upshape)}.  Then for every $\e>0$ there is a proper closed
subvariety $Z=Z_\e(E)\subset E^2$ so that 
\[
  \gcd\bigl(D_P,D_Q\bigr) \le  \bigl(H(P)\cdot H(Q)\bigr)^\e
  \qquad\text{for all $(P,Q)\in E^2(\QQ)\setminus Z$.}
\]
\par
The exceptional set~$Z$ consists of a finite number of translates of
proper algebraic subgroups of~$E^2$. If~$E$ does not have~CM, then we
can say more precisely~$Z$ is a finite union of translates of the
subgroups
\[
  \bigl\{(mT,nT)\in E^2 : T\in E\bigr\}
  \quad\text{with $(m,n)\in\ZZ^2$ satisfying $m^2+n^2\le\frac{1}{2\e}$.}
\]
{\upshape(}A similar statement holds if~$E$ has~CM, with~$m$ and~$n$
replaced by more general isogenies.{\upshape)}
\end{theorem}

\begin{example}
Let $E/\QQ$ be an elliptic curve and $P\in E(\QQ)$ a point of infinite
order. With notation as in Theorem~\ref{theorem:ECoverQ}, the
\emph{elliptic divisibility sequence} (EDS) associated to~$P$ is the
sequence of integers~$(D_{nP})_{n\ge1}$.  (For further information
about elliptic divisibility sequences, including a
not-quite-equivalent alternative definition, see
\cite{EGW,EPSW,EW1,McKinnon,SHIP,JHSECFF,JHSEDS,SS,SW,Ward1,Ward2}.)
These sequences have the property that if $m|n$, then $D_{mP}|D_{nP}$,
whence their name. Now let~$P$ and~$Q$ be independent points in~$E(\QQ)$.
Then Theorem~\ref{theorem:ECoverQ} 
implies that there is a constant $C=C_\e(E,P,Q)$ so that
\[
  \gcd\bigl(D_{mP},D_{nQ}\bigr) \le C\max\bigl\{D_{mP},D_{nQ}\bigr\}^\e
  \qquad\text{for all $m,n\ge1$.}
\]
(Note that since~$P$ and~$Q$ are independent, there are only finitely 
many multiples~$(mP,nQ)$ that lie on any fixed curve in~$E^2$.
We are also using Siegel's theorem~\cite[IX.3.3]{AEC},
which says that $2\log D_{nP}\sim h(nP)$ as $n\to\infty$.)
\end{example}

Our final example is the amusing observation that Vojta's conjecture
allows us to mix greatest common divisors of a multiplicative group
with an elliptic curve. The following result, although far from the
most general, gives a flavor of what is can be proven.
Again, an unconditional proof would be quite interesting.

\begin{theorem}
\label{theorem:mixedEandGm}
Let $E/\QQ$ be an elliptic curve and let~$S$ be a finite set of
rational primes.  Assume that Vojta's conjecture is true for
\text{$E\times\PP^1$} blown up at~$(O,1)$.  Then for every $\e>0$
there is a constant $C=C(E,S,\e)$ so that
\[
  \gcd(D_{Q},b-1) \le C\cdot\max\{D_{Q},b\}^\e \quad
  \text{for all $Q\in E(\QQ)$ and $b\in\ZZ_S^*\cap\ZZ$.}
\]
{\upshape(}By convention, we define the greatest common divisor of two 
rational numbers to be the  greatest common divisor of their
numerators.{\upshape)} In particular, if $P\in E(\QQ)$ is a point
of infinite order and if $a\ge2$ is an integer, then
\[
  \gcd(D_{nP},a^m-1) \le \max\{D_{nP},a^m\}^\e
\]
provided that $\max\{m,n\}$ is sufficiently large.
\end{theorem}

\section{Generalized gcds and blowups}
\label{section:generalizedgcddef}

We set the following notation, which will remain fixed throughout this
paper.  For definitions and normalizations related to absolute values
and heights, see~\cite[Part~B]{HS} or \cite[Chapters~2,~3]{LangDG}.

\begin{notation}
\item[$k$]
a number field.
\item[$M_k$]
a complete set of absolute values on~$k$. For $v\in M_k$, we define
$v^+(\a)=\max\{v(\a),0\}$, and we assume that the absolute values are
normalized so that $h(\a)=\sum_{v\in M_k} v^+(\a)$ is the absolute
logarithmic Weil height of~$\a$.  We denote by~$M_k^0$, respectively
by~$M_k^\infty$, the set of nonarchimedean, respectively archimedean,
places in~$M_k$.
\item[$S$]
a finite set of places of~$k$, including all of the archimedean places.
\item[$X/k$]
a smooth projective variety defined over $k$.
\item[$h_{X,D}$]
an absolute logarithmic Weil height on~$X$ with respect to
the divisor~$D$. 
\item[$\l_{X,D}$]
an absolute logarithmic local height on~$X$ with respect to
the divisor~$D$. 
\end{notation}

Let $a,b\in\ZZ$. The greatest common divisor of~$a$ and~$b$ is given
by the formula
\begin{align*}
  \log\gcd(a,b) 
  &= \sum_{p~\text{prime}}
       \min\bigl\{\ord_p(a),\ord_p(b)\bigr\}\log p\\
  &= \sum_{v\in M_\QQ^0}  \min\bigl\{v(\a),v(\b)\bigr\}.
\end{align*}
If~$a$ and~$b$ are rational numbers, rather than integers, then we can
compute the gcd of their numerators by using~$v^+$ in place of~$v$,
and having done this, there is no reason to restrict ourselves to the
nonarchimedean places. Moving from~$\QQ$ to the number field~$k$, we
follow~\cite{CZ2} and define the \emph{generalized}
(\emph{logarithmic}) \emph{greatest common divisor} of~$\a,\b\in k$ to
be the quantity
\[
  h_{\gcd}(\a,\b) = \sum_{v\in M_k} \min\bigl\{v^+(\a),v^+(\b)\bigr\}.
\]
In particular, if~$\a,\b\in\ZZ$, then $h_{\gcd}(\a,\b)=\log\gcd(\a,\b)$.

A fancier way to view the function
\[
  v^+ : k \longrightarrow [0,\infty]
\]
is as the local height function on~$\PP^1(k)$ with respect to the
divisor~$(0)$, where we identify $k\cup\{\infty\}$ with~$\PP^1(k)$ and
set $v^+(\infty)=0$. We would like to find  a similar 
height theoretic interpretation for the function
\[
  G:\PP^1(k)\times\PP^1(k) \longrightarrow [0,\infty],\qquad
  (\a,\b) \longmapsto \min\bigl\{v^+(\a),v^+(\b)\bigr\},
\]
that appears in the definition of the generalized greater common
divisor.  Intuitively,~$G(\a,\b)$ is large if and only if the
point~$(\a,\b)$ is $v$-adically close to the point~$(0,0)$. This 
resembles the intuitive characterization of a local height function,
\[
  \l_{X,D}(P,v) = -\log(\text{$v$-adic distance from $P$ to $D$}),
\]
except that~$(0,0)$ is not a divisor on~$(\PP^1)^2$. However, there is a
general theory that associates a local height function~$\l_{X,Y}(P,v)$
to any subvariety~$Y$ of~$X$, or more generally to any closed
subscheme~$Y$, see~\cite{JHSADF} or~\cite[\S5]{Vojta}. For our
purposes, it is convenient to use an equivalent formulation in
terms of blowups.
\par
Continuing with our example, let $X=(\PP^1)^2$, let $\pi:\Xtilde\to
X$ be the blowup of~$X$ at the point~$(0,0)$, and let $E=\pi^{-1}(0,0)$
be the exceptional divisor of the blowup.  Then it is an easy exercise
using explicit equations (or see~\cite[Lemma~2.5.2]{Vojta}) to verify
that a local height function on~$\Xtilde$ for the divisor~$E$ is given
by the formula
\begin{multline*}
  \l_{\Xtilde,E}(\pi^{-1}(\a,\b),v)
  = \min\bigl\{v^+(\a),v^+(\b)\bigr\} \\
  \quad\text{for all $(\a,\b)\in X(k)\setminus(0,0)$.}
\end{multline*}
Adding these local heights gives the global formula
\[
  h_{\gcd}(\a,\b)
  = \sum_{v\in M_k}   \l_{\Xtilde,E}(\pi^{-1}(\a,\b),v)
  = h_{\Xtilde,E}(\pi^{-1}(\a,\b)).
\]
In other words, the (generalized) logarithmic gcd of~$\a$ and~$\b$ is
equal to the Weil height of~$(\a,\b)$ on a blowup of~$(\PP^1)^2$
with respect to the exceptional divisor of the blowup. This identification
allows us to bring the machinery of heights to bear on problems
concerning greatest common divisors, and in particular allows us to
apply Vojta's conjecture to such problems.
\par
Having identified~$h_{\gcd}(\a,\b)$ with the Weil height on a
particular blowup, it is natural to generalize the notion of
greatest common divisor to arbitrary varieties blown up up along
arbitrary subvarieties.

\begin{definition}
Let $X/k$ be a smooth variety and let $Y/k\subset X/k$ be a 
subvariety of codimension $r\ge2$. Let $\pi:\Xtilde\to X$ be the
blowup of~$X$ along~$Y$, and let $\Ytilde=\pi^{-1}(Y)$ be the
exceptional divisor of the blowup. For $P\in X\setminus Y$, we
let~$\Ptilde=\pi^{-1}(P)\in\Xtilde$.
\par
The \emph{generalized} (\emph{logarithmic}) \emph{greatest common
divisor of the point \text{$P\in(X\setminus Y)(k)$} with respect
to~$Y$} is the quantity
\[
  h_{\gcd}(P;Y) = h_{\Xtilde,\Ytilde}(\Ptilde).
\]
\end{definition}

\begin{example}
\label{example:gcdPn000}
Let $X=\PP^n$ and let $Y=[1,0,0,\ldots,0]$. For $\bfx\in\PP^n(\QQ)$,
choose homogeneous coordinates $\bfx=[x_0,x_1,\ldots,x_n]$
with $x_i\in\ZZ$ and $\gcd(x_0,\ldots,x_n)=1$. Then
\begin{multline*}
  h_{\gcd}(\bfx;Y) = \log \gcd(x_1,x_2,\ldots,x_n)+O(1) \\
  \qquad\text{for $\bfx=[x_0,x_1,\ldots,x_n]\in\PP^n(\QQ)$.}
\end{multline*}
\end{example}

\begin{example}
\label{example:gcdPnY}
Again let $X=\PP^n$ and let~$Y$ be a subvariety of codimension $r\ge2$
defined by the vanishing of a collection of homogeneous polynomials
$f_1,f_2,\ldots,f_t\in\ZZ[X_0,\ldots,X_n]$.   Then for all points
$\bfx=[x_0,x_1,\ldots,x_n]\in\PP^n(\QQ)$ written with normalized
homogeneous coordinates as in Example~\ref{example:gcdPn000}, we have
\[
  h_{\gcd}(\bfx;Y) = \log \gcd\bigl(f_1(\bfx),\ldots,f_t(\bfx)\bigr)+O(1).
\]
Compare the righthand side of this formula with the lefthand side
of condition~(2) in Theorem~\ref{theorem:PnoverQ}. This will allow us
to reformulate  Theorem~\ref{theorem:PnoverQ} in terms of heights on
blown up varieties and thence to apply Vojta's conjecture.
\end{example}

\begin{example}
\label{example:gcdE2}
Let $E/\QQ$ be an elliptic curve given by a (minimal) Weierstrass
equation, let $X=E^2$, let $Y=\{(O,O)\}$, and let
$\pi_1,\pi_2:X\to E$ denote the two projections.  The 
square of the ideal sheaf~$\Ical_Y$
of~$Y$ is generated locally by the two functions~$\pi_1^*(x^{-1})$
and~$\pi_2^*(x^{-1})$,
\[
  \Ical_Y^2 = \pi_1^*(x^{-1})\Ocal_{X,Y} + \pi_2^*(x^{-1})\Ocal_{X,Y}.
\]
Hence the greatest common divisor of a  point~$(P,Q)\in X(\QQ)$
with respect to $Y=\{(O,O)\}$ is given by
\begin{align}
  \label{equation:hgcd=gcdDPDQonE2}
  h_{\gcd}\bigl((P,Q);Y\bigr)
  &= \sum_{v\in M_\QQ} \frac{1}{2}\min\bigl\{v^+(x_P^{-1}),v^+(x_Q^{-1})\}
    \notag\\
  &= \log\gcd(D_P,D_Q),
\end{align}
where recall (cf{.}~Theorem~\ref{theorem:ECoverQ})
that for $P\in E(\QQ)$, we write $x_P=A_P/D_P^2$.
\end{example}

\section{Vojta's conjecture}
\label{section:vojtaconjecture}
We recall the statement of Vojta's 
conjecture~\cite[Conjecture~3.4.3]{Vojta}.

\begin{conjecture}[Vojta~\cite{Vojta}]
\label{conjecture:vojta}
Set the following notation:
\begin{notation}
\item[$k$]
a number field.
\item[$S$]
a finite set of places of~$k$.
\item[$X/k$]
 a smooth projective variety.
\item[$A$]
an ample divisor on~$X$.
\item[$K_X$]
 a canonical divisor on~$X$.
\end{notation}
Then for every $\e>0$ there exists a proper Zariski-closed subset
$Z=Z(\e,X,A,D,k,S)$ of~$X$ and a constant $C_\e=C_\e(X,A,D,k,S)$ so
that
\begin{multline}
  \label{equation:vojtaconj}
  \smash{\sum_{v\in S}}\ \l_{X,D}(P,v) + h_{X,K_X}(P) \le \e h_{X,A}(P)+C_\e \\
  \qquad\text{for all $P\in X(k)\setminus Z$.}
\end{multline}  
\end{conjecture}

\begin{remark}
Vojta's conjecture contains the additional statement
that aside from a set of dimension zero, the set~$Z$ may be chosen
independently of the field~$k$ and the set of places~$S$. In other
words, there is a set~$Z_0=Z_0(\e,X,A,D)$ so that for any finite
extension~$k'/k$ and any finite set of places~$S'$ of~$k'$, there is a
finite set of \emph{points} $Z_1=Z_1(\e,X, A,D,k',S')$ so
that~\eqref{equation:vojtaconj} holds for all $P\in X(k')$ with
$P\notin Z=Z_0\cup Z_1$. We will be working over a single
number field, so we will not need this stronger version.
\end{remark}

\begin{remark}
In Vojta's conjecture and throughout this paper, when we say that a constant
depends on a divisor~$D$ on a variety~$X$, we assume that both 
global and local heights~$h_{X,D}$ and~$\l_{X,D}$ have been chosen and that
the constant in question may depend on this choice. 
\end{remark}

\begin{definition}
With notation as in the statement of Conjecture~\ref{conjecture:vojta},
we let
\begin{align*}
  h_{X,D,S}(P) &= \sum_{v\in S} \l_{X,D}(P,v), \\
  h_{X,D,S}'(P) &= \sum_{v\notin S} \l_{X,D}(P,v).
\end{align*}
This corresponds to Vojta's notation~\cite{Vojta} via
$m_S(D,P)=h_{X,D,S}(P)$ and $N_S(D,P)=h_{X,D,S}'(P)$.  Making an
analogy with Nevanlinna theory, Vojta calls~$m_S(D,P)$ the ``proximity
function'' and~$N_S(D,P)$ the ``counting function.''  Then Vojta's
fundamental inequality~\eqref{equation:vojtaconj} becomes the succinct
statement
\begin{equation}
  \label{equation:vojtaconj1}
  h_{X,D,S}(P) + h_{X,K_X}(P) \le \e h_{X,A}(P)+C_\e
  \qquad\text{for all $P\in X(k)\setminus Z$.}
\end{equation}  
\end{definition}

\section{Applying Vojta's conjecture to blowups}
\label{section:applyvojtatoblowups}

Let $X/k$ be a smooth variety and let $Y/k\subset X/k$ be a smooth
subvariety of codimension $r\ge2$. Let $\pi:\Xtilde\to X$ be the
blowup of~$X$ along~$Y$, and let $\Ytilde=\pi^{-1}(Y)$ be the
exceptional divisor of the blowup. For $P\in X\setminus Y$, we
let~$\Ptilde=\pi^{-1}(P)\in\Xtilde$. A nice property of blowups of
smooth varieties along smooth subvarieties is that it is easy to
describe the canonical bundle on the
blowup~\cite[Exercise~II.8.5]{Hartshorne},
\[
  K_{\Xtilde} \sim \pi^*K_X + (r-1)\Ytilde.
\]
(Here $\sim$ denotes linear equivalence.)  We also observe that if~$A$
is an ample divisor on~$X$, then there exists an integer~$N$ so that
\text{$-\Ytilde+N\pi^*A$} is ample on~$\Xtilde$.  This follows from
the Nakai-Moishezon Criterion~\cite[Theorem~A.5.1]{Hartshorne}.  We
choose such an~$N$ and let
\[
  \Atilde=-\frac{1}{N}\Ytilde+\pi^*A\in\Div(\Xtilde)\otimes\QQ,
\]
so~$\Atilde$ is in the ample cone of~$\Xtilde$.
\par
We make the following assumption:
\begin{equation}
  \label{equation:fanoassumption}
  \left(
   \parbox{.65\hsize}{\noindent
      The anticanonical divisor $-K_X$ is 
      a normal crossings divisor and $\Support(K_X)\cap Y=\emptyset$.}
   \right)
\end{equation}
(In practice, it suffices to assume that some multiple of $-K_X$ is a
normal crossings divisor. The case $K_X=0$ is also permitted.)  With
notation as above and under the
assumption~\eqref{equation:fanoassumption}, we apply Vojta's
conjecture to the variety~$\Xtilde$ and the divisor $D=-\pi^*K_X$
to obtain the inequality
\[
  h_{\Xtilde,-\pi^*K_X,S}(\Ptilde)
  + h_{\Xtilde,K_\Xtilde}(\Ptilde)
  \le \e h_{\Xtilde,\Atilde}(\Ptilde) + C_\e
  \qquad\text{for all $\Ptilde\in \Xtilde(k)\setminus \Ztilde$.}
\]
Substituting
\[
  K_\Xtilde=\pi^*K_X + (r-1)\Ytilde
  \qquad\text{and}\qquad
  \Atilde=-\frac{1}{N}\Ytilde+\pi^*A
\]
and using functorial properties of height functions, we obtain
\begin{multline*}
  -h_{X,K_X,S}(P) + h_{X,K_X}(P) + (r-1)h_{\Xtilde,\Ytilde}(\Ptilde) \\
  \le \e h_{X,A}(P) - \frac{\e}{N}h_{\Xtilde,\Ytilde}(\Ptilde) + C_\e 
  \qquad\text{for all $P\in X(k)\setminus Z$,}
\end{multline*}
where we have written $Z=\pi(\Ztilde)$. The two leftmost terms may be
combined using $h_{X,D,S}+h_{X,D,S}'=h_{X,D}$, which yields
\begin{multline*}
  h_{X,K_X,S}'(P) + \left(r-1+\frac{\e}{N}\right)h_{\Xtilde,\Ytilde}(\Ptilde)
  \le \e h_{X,A}(P) + C_\e \\
  \qquad\text{for all $P\in X(k)\setminus Z$.}
\end{multline*}
Finally, a small amount of algebra,
the definition $h_{\gcd}(P;Y)=h_{\Xtilde,\Ytilde}(\Ptilde)$,
and setting $\d=\e/N$ gives the
following result, where for the convenience of the reader we restate
all of our assumptions.

\begin{theorem}
\label{theorem:vojta=>hgcdbound}
Let $X/k$ be a smooth variety, let~$A$ be an ample divisor on~$X$, and
let $Y/k\subset X/k$ be a smooth subvariety of codimension
$r\ge2$. Assume that~$-K_X$ is a normal crossings divisor whose
support does not intersect~$Y$.  Assume further that Vojta's
conjecture is true {\upshape(}at least for the blowup $\pi:\Xtilde\to
X$ of~$X$ along~$Y$ and for the divisor $D=-\pi^*K_X${\upshape)}.  Then
for every finite set of places~$S$ and every $0<\e<r-1$ there is a proper
closed subvariety $Z=Z(\e,X,Y,A,k,S)\subsetneq X$, a constant
$C_\e=C_\e(X,Y,A,k,S)$, and a constant $\d=\d(X,Y,A)$ so that
\begin{multline}
  \label{equation:generalhgcdbound}
  h_{\gcd}(P;Y) \le \e h_{X,A}(P) + 
      \frac{1}{r-1+\d\e}h_{X,-K_X,S}'(P) + C_\e\\
  \qquad\text{for all $P\in X(k)\setminus Z$.}
\end{multline}
\end{theorem}

\section{Proofs of Theorems~\ref{theorem:PnoverQ},~\ref{theorem:ECoverQ},
and~\ref{theorem:mixedEandGm}}
\label{section:proofsof3theorems}

In this section we show how our main result
(Theorem~\ref{theorem:vojta=>hgcdbound}) can be used to prove the three
special cases stated in Section~\ref{section:examplesoverQ}.

\begin{proof}[Proof of Theorem~\ref{theorem:PnoverQ}]
We apply  Theorem~\ref{theorem:vojta=>hgcdbound} to the following data:
\begin{align*}
  X & = \PP^n, \\
  Y &= \{f_1=f_2=\cdots=f_t=0\} \subset \PP^n, \\
  K_X &= -\sum_{i=0}^n H_i,
    \qquad\text{where $H_i=\{X_i=0\}\in\Div(\PP^n)$,}\\
  A &= H_0.
\end{align*}
Notice that~$-K_X$ is a normal crossings divisor and 
that~$Y$ is disjoint from the support of~$-K_X$ by assumption.
For $P\in\PP^n(\QQ)$, let $\bfx=(x_0,\ldots,x_n)\in\ZZ^{n+1}$
with $\gcd(x_i)=1$ be  normalized homogeneous coordinates for~$P$. Then
by definition of the Weil height we have
\begin{equation}
  \label{equation:ht1}
  h_{X,A} = \log \max\bigl\{|x_0|,\ldots,|x_n|\bigr\},
\end{equation}
and Example~\ref{example:gcdPnY} says that
\begin{equation}
  \label{equation:ht2}
  h_{\gcd}(P;Y) = \log \gcd\bigl(f_1(\bfx),\ldots,f_t(\bfx)\bigr).
\end{equation}
(All height equalities up to~$O(1)$.)
Further, by definition of the $S$-part of the height, we have
\[
  h_{X,H_i,S}'(P) = \sum_{v\notin S} v^+(x_i) = \log|x_i|_S',
\]
so
\begin{equation}
  \label{equation:ht3}
  h_{X,-K_X,S}'(P) = \sum_{i=0}^n h_{X,H_i,S}'(P)
      = \log|x_0x_1\cdots x_n|_S'.
\end{equation}
\par
We now substitute~\eqref{equation:ht1},~\eqref{equation:ht2}
and~\eqref{equation:ht3} into the
inequality~\eqref{equation:generalhgcdbound} of
Theorem~\ref{theorem:vojta=>hgcdbound} to obtain
\begin{multline*}
  \log \gcd\bigl(f_1(\bfx),\ldots,f_t(\bfx)\bigr) \\
  \le \e \log \max\bigl\{|x_0|,\ldots,|x_n|\bigr\} 
  + \frac{1}{r-1+\d\e}\log|x_0x_1\cdots x_n|_S'
  + C_\e \\
  \text{for all $P=[\bfx]\in\PP^n(\QQ)\setminus Z$.}
\end{multline*}
Exponentiating this inequality completes the proof of
Theorem~\ref{theorem:PnoverQ}, once we observe that the exceptional
set~$Z$ is contained in some hypersurface, so may be replaced by the
zero set of a single nonzero polynomial.
\end{proof}

\begin{proof}[Proof of Theorem~\ref{theorem:ECoverQ}]
Let $\pi_1,\pi_2:E\times E\to E$ be the two projections. 
We apply  Theorem~\ref{theorem:vojta=>hgcdbound} to the following data:
\[
  X  = E\times E,\qquad Y=\{(O,O)\},\qquad K_X=0,\qquad
  A=\pi_1^*(O)+\pi_2^*(O).
\]
We compute
\begin{alignat}{2}
  \label{equation:ht4}
  h_{X,A}(P,Q) &= h_{E\times E,\pi_1^*(O)+\pi_2^*(O)}(P,Q)
    &\qquad&\text{definition of $X$ and $A$,} \notag\\
  &= h_{E,O}(P)+h_{E,O}(Q) + O(1)
    &&\text{functoriality of heights.}
\end{alignat}
Next we recall from~\eqref{equation:hgcd=gcdDPDQonE2} in
Example~\ref{example:gcdE2} that the generalized greatest common
divisor of~$(P,Q)$ with respect to~$(O,O)$ is given by 
\begin{equation}
  \label{equation:ht5}
  h_{\gcd}\bigl((P,Q);(O,O)\bigr)
  = \log\gcd(D_P,D_Q).
\end{equation}
Substituting~\eqref{equation:ht4} and~\eqref{equation:ht5}
into inequality~\eqref{equation:generalhgcdbound} of
Theorem~\ref{theorem:vojta=>hgcdbound} yields (note~$K_X=0$,
so the~$h'_{X,-K_X,S}$ term disappears)
\begin{multline*}
  \log\gcd(D_P,D_Q) \le \e \bigl(h_{E,O}(P)+h_{E,O}(Q)) + C_\e \\
  \qquad\text{for all $(P,Q)\in E^2(\QQ)\setminus Z$.}
\end{multline*}
Exponentiating gives the first part of Theorem~\ref{theorem:ECoverQ}.
\par
It remains to describe the exceptional set~$Z$. Let~$\G\subset Z$ be
an irreducible component of~$Z$ such that
\begin{multline}
  \label{equation:ht6}
  \log\gcd(D_P,D_Q) \ge \e \bigl(h(P)+h(Q)\bigr) + C_\e \\
  \qquad\text{for infinitely many $(P,Q)\in\G(\QQ)$,}
\end{multline}
where to ease notation we let~$h(P)=h_{E,O}(P)$.
Faltings' theorem~\cite{Faltings} tells us that~$\G$ is a translate
of an abelian subvariety of~$E^2$, i.e.~$\G$ is a translate
of an elliptic curve. If~$E$ does not have~CM, then the abelian subvarieties
of~$E^2$ are precisely the curves
\[
  \G_{n_1,n_2} = \bigl\{(n_1T,n_2T):T\in E\bigr\}
  \qquad\text{for $n_1,n_2\ge0$ with $\gcd(n_1,n_2)=1$.}
\]
Thus the assumption that~$\G$ contains infinitely many points
satisfying~\eqref{equation:ht6} implies that there is a fixed pair of
integers~$(n_1,n_2)$ as above and a fixed pair of points~$(R_1,R_2)\in
E^2(\QQ)$ so that 
\[
  \G=\G_{n_1,n_2}+(R_1,R_2)
  =\bigl\{(n_1T+R_1,n_2T+R_2) : T\in E\bigr\}.
\]
Hence
\begin{align}
  \label{equation:ht7}
  \log\gcd(D_{n_1T+R_1},D_{n_2T+R_2}) 
    &\ge \e\bigl(h(n_1T+R_1)+h(n_2T+R_2)\bigr) + O(1) \notag\\
    &= \e(n_1^2+n_2^2)h(T) + O\left(\sqrt{h(T)}\,\right) \notag\\
    &\omit\hfil\text{for infinitely many $T\in E(\QQ)$.}\qquad
\end{align}
Here the big-$O$ constant may depend on~$(R_1,R_2)$ and
on~$(n_1,n_2)$, as long as it is independent of~$T$. We have also used
the positivity and quadratic nature of the height (\cite[VIII~\S9]{AEC})
in the form
\[
  h(nT+R)=n^2h(T)+O_{E,R}\left(\sqrt{h(T)}\,\right).
\]
\par
It remains to bound $\gcd(D_{n_1T+R_1},D_{n_2T+R_2})$.  Since
$\gcd(n_1,n_2)=1$ by assumption, we can choose integers~$(u_1,u_2)$ with
\text{$u_1n_1+u_2n_2=1$} and set \text{$R_3=u_1R_1+u_2R_2$}. Note
that~$R_3$ is independent of~$T$.  Let~$p$ be a prime. Working
in~$E(\QQ_p)$, we have
\begin{align*}
  p^e|\gcd(&D_{n_1T+R_1},D_{n_2T+R_2})\\
   {}\Longleftrightarrow{}
  &n_1T+R_1\equiv O\pmod{p^e}\quad\text{and}\quad n_2T+R_2\equiv O\pmod{p^e} \\
  {}\Longrightarrow{}
  &T+R_3= u_1(n_1T+R_1) + u_2(n_2T+R_2) \equiv O \pmod{p^e} \\
  {}\Longrightarrow{}
  &p^e|D_{T+R_3}.
\end{align*}
Thus $\gcd(D_{n_1T+R_1},D_{n_2T+R_2})$ divides~$D_{T+R_3}$, so
\begin{align}
  \label{equation:gcdD1D2hT}
  \log\gcd(D_{n_1T+R_1},D_{n_2T+R_2})
  &\le \log D_{T+R_3} \notag\\
  &\le h(T+R_3) \notag\\
  &\le h(T) + O(\sqrt{h(T)}\,).
\end{align}
Combining~\eqref{equation:ht7} and~\eqref{equation:gcdD1D2hT} yields
\[
  h(T) \ge 
  \e(n_1^2+n_2^2)h(T) + O\bigl(\sqrt{h(T)}\,\bigr) \quad
  \text{for infinitely many $T\in E(\QQ)$.}
\]
Letting $h(T)\to\infty$, we conclude that
\begin{equation}
  \label{equation:1geen1n2}
  1 \ge \e(n_1^2+n_2^2).
\end{equation}
This completes the proof of Theorem~\ref{theorem:ECoverQ}
once we observe that the height function~$H(P)$ used in the
statement of Theorem~\ref{theorem:ECoverQ}
satisfies $\log H(P) = 2h_{E,O}(P)$.
\end{proof}

\begin{proof}[Proof of Theorem~\ref{theorem:mixedEandGm}]
This time we apply Theorem~\ref{theorem:vojta=>hgcdbound} with
\begin{align*}
  X &= E\times\PP^1,\\
  A &= \pi_1^*(O) + \pi_2^*(\infty),\\
  K_X &= -\pi_2^*(0)-\pi_2^*(\infty),\\
  Y &= \{(O,1)\},
\end{align*}
where $\pi_1:X\to E$ and $\pi_2:X\to\PP^1$ are the projections.
Then for any $(Q,b)\in E(\QQ)\times\ZZ$ we have  
\begin{align*}
  h_{X,A}(Q,b) &= h_{E,O}(Q) + h(b)\\
  h_{\gcd}\bigl((Q,b);(0,1)\bigr) &= \log\gcd(D_Q,b-1).
\end{align*}
Further, if~$b\in\ZZ_S^*$, then
\[
  h_{X,-K_X,S}'(Q,b) = h_{\PP^1,(0),S}'(b) + h_{\PP^1,(\infty),S}'(b) = 0.
\]
Thus Theorem~\ref{theorem:vojta=>hgcdbound} yields
\begin{multline*}
  \log\gcd(D_Q,b-1) \le \e\bigl(h_{E,O}(Q) + h(b)\bigr) + O(1)\\
  \text{for $(Q,b)\in E(\QQ)\times\ZZ_S^*$ with $(Q,b)\notin Z$.}
\end{multline*}
Siegel's theorem~\cite[IX.3.3]{AEC} says that $h_{E,O}(Q)\sim\log D_Q$
as $h_{E,O}(Q)\to\infty$, so exponentiating and adjusting~$\e$ gives
\begin{multline*}
  \gcd(D_Q,b-1) \le C\cdot \max\bigl(D_Q,b\bigr)^\e \\
  \text{for $(Q,b)\in E(\QQ)\times\ZZ_S^*$ with $(Q,b)\notin Z$.}
\end{multline*}
It remains to deal with the exceptional set~$Z$. It suffices to
consider an irreducible component $\G\subset Z$ of dimension~1 with
\begin{multline}
  \label{equation:ht8}
  \log\gcd(D_Q,b-1) \ge \e\bigl(h_{E,O}(Q) + h(b)\bigr) + O(1)\\
  \text{for infinitely many $(Q,b)\in \bigl(E(\QQ)\times\ZZ_S^*\bigr)\cap \G$.}
\end{multline}
In particular,
$\#\G(\QQ)=\infty$, so Faltings' theorem~\cite{FaltingsMC} reduces us
to the case that~$\G$ has genus~0 or~1.  If either~$\pi_1(\G)$
or~$\pi_1(\G)$ consists of a single point, it suffices to adjust the
constant, so we assume that $\pi_1(\G)=E$ and $\pi_2(\G)=\PP^1$. In particular,
the fact that $\pi_1(\G)=E$ implies that~$\G$ cannot have genus~0, so
we are reduced to the case that~$\G$ has genus~$1$.
\par
The fact that~$\G$ satisfies~\eqref{equation:ht8} implies
that $\pi_2(\G)\cap\ZZ_S^*$ is infinite. In other words, the
map
\[
  \pi_2:\G(\QQ)\longrightarrow \QQ\cup\{\infty\}
\]
takes on infinitely many $S$-unit values. But~$\G(\QQ)$ is the Mordell-Weil
group of an elliptic curve, so Siegel's theorem~\cite[IX.3.2.2]{AEC}
says that this is not possible  (indeed, it is not even possible to
take on infinitely many $S$-integral values). This completes the proof 
that the exceptional set may be taken to be a finite set of points, and
hence may be eliminated entirely by adjusting the constants.
\end{proof}

\section{Divisibility sequences and algebraic groups}
\label{section:divseqsandalggps}

A \emph{divisibility sequence} is a sequence of integers~$(a_n)_{n\ge1}$
with the property that
\[
  m|n\Longrightarrow a_m|a_n.
\]
We have already briefly discussed the divisibility
sequences~$(D_{nP})$ associated to a point of infinite order~$P$ on an
elliptic curve~$E(\QQ)$. Other familiar divisibility sequences include
sequences of the form \text{$(a^n-b^n)$} and the Fibonacci
sequence~$(F_n)$. There are many natural ways to generalize the notion
of divisibility sequence, for example by replacing divisbility of
positive integers with divisibility of ideals in a ring. In the most
abstract formulation, one might define a divisibility sequence as
simply an order-preserving maps between two partially ordered sets
(posets). In this section we  restrict our attention to classical
divisibility sequences of rational integers, but the reader should be
aware that virtually everything that we say can be easily generalized
(albeit at the cost of some notational inconvenience) to the partially
ordered set of integral ideals in number fields, and in some cases
to other Dedekind domains or even more general rings.
\par
The divisibility sequence \text{$(a^n-b^n)_{n\ge1}$} is naturally
associated to the rank one subgroup of~$\GG_m(\QQ)$ generated
by~$a/b$, just as the divisibility sequence~$(D_{nP})_{n\ge1}$ comes from the
rank one subgroup of~$E(\QQ)$ generated by~$P$. This suggests creating
divisibility sequences from other algebraic groups~$G$ defined
over~$\QQ$. In order to make this precise, we need to choose a model
over~$\ZZ$, although a  a different choice of model only changes the
sequence at finitely many primes. 

\begin{definition}
Let $\Gcal/\ZZ$ be a group scheme over~$\ZZ$, let
$\Ocal\subset\Gcal(\ZZ)$ be the identity element of~$\Gcal$, and let
$\Pcal\in\Gcal(\ZZ)$ be a nonzero section. We associate to~$\Pcal$
a positive integer~$D_\Pcal$ by the condition
\[
  \ord_p(D_{\Pcal}) = (\Pcal\cdot\Ocal)_p
  \qquad\text{for all primes $p$,}
\]
where in general $(\Pcal_1\cdot\Pcal_2)_p$ denotes the arithmetic intersection
index of the sections~$\Pcal_1$ and~$\Pcal_2$ on the fiber over~$p$.
\par
Equivalently, 
let $\Ical_\Ocal$ be the ideal sheaf of~$\Ocal\subset\Gcal$, where
we identify the section~$\Ocal$ with its image~$\Ocal(\ZZ)$, taken with
the induced reduced subscheme structure. Then $\Pcal^*(\Ical_\Ocal)$
is an ideal sheaf on~$\Spec(\ZZ)$, i.e. it is an ideal of~$\ZZ$.
Then~$D_\Pcal$ is determined by the condition that it generates this ideal,
\[
  D_{\Pcal}\cdot\ZZ = (\Pcal)^*(\Ical_\Ocal).
\]
\end{definition}

These~$D_\Pcal$ values are closely associated to certain generalized
greatest common divisors.

\begin{proposition}
\label{proposition:logDPlehgcd}
Let~$\Gcal/\ZZ$ be a group scheme, let~$G=\Gcal\times_\ZZ\QQ$ be the
associated algebraic group over~$\QQ$, and let $\r:\Gcal(\ZZ)\to
G(\QQ)$ denote restriction to the generic fiber,and let
$O=\r(\Ocal)\in G(\QQ)$ be the identity element of~$G$. Then
\[
  \log D_{\Pcal}  \le  h_{\gcd}(\r(\Pcal);O) + O(1)
  \qquad\text{for all $\Pcal\in\Gcal(\ZZ)$.}
\]
{\upshape(}In principle, the height function might depend on the
choice of a completion and projective embedding
of~$G$. However, these only affect~$h_{\gcd}(\,\cdot\,;O)$ 
up to~$O(1)$.{\upshape)}
\end{proposition}
\begin{proof}
This is just a matter of unsorting the definitions and
decomposing~$h_{\gcd}$ into a sum of local heights. With the obvious
notation, we find that
\begin{multline*}
  \l_{\gcd}(\r(\Pcal);O;v) 
  = \l_{\tilde G,\tilde O}(\r(\Pcal),v)
  = v(D_{\Pcal})\\
  \text{for all nonarchimedean places~$v$.}
\end{multline*}
This gives the stated result, with the contributions
from the (nonnegative) archimedean local heights giving an inequality,
rather than an equality.
\end{proof}

We next show that a sequence of the form $(D_{n\Pcal})_{n\ge1}$ is
a divisibility sequence.

\begin{proposition}
\label{proposition:Dnisdivseq}
Let $\Gcal/\ZZ$ be a group scheme and let~$\Pcal\in\Gcal(\ZZ)$ be a
point {\upshape(}section{\upshape)} of infinite order.  Then the
sequence~$(D_{n\Pcal})_{n\ge1}$ 
is a divisibility sequence.  We call it the \emph{divisibility
sequence associated to~$\Pcal$ {\upshape(}and~$\Gcal${\upshape)}.}
\end{proposition}
\begin{proof}
For each integer $n\ge1$, let $\m_n:\Gcal\to\Gcal$ be the
$n^{\text{th}}$-power morphism. The section~$n\Pcal\in\Gcal(\ZZ)$ is
the composition
\[
  \begin{CD}
  \Spec(\ZZ) @>\Pcal>> \Gcal @>\m_n>> \Gcal.
  \end{CD}
\]
Now let $m|n$, say $n=mr$. Then
\begin{alignat*}{2}
  D_{n\Pcal}\cdot\ZZ  &= (n\Pcal)^*(\Ical_\Ocal)
    &\qquad&\text{by definition of $D_{n\Pcal}$,}\\
  &=(\m_n\circ\Pcal)^*(\Ical_\Ocal)
    &&\text{since $n\Pcal=\m_n\circ\Pcal$ as maps,}\\
  &=(\m_r\circ\m_m\circ\Pcal)^*(\Ical_\Ocal)
    &&\text{since $\m_n=\m_{rm}=\m_r\circ\m_m$,}\\
  &=(\m_m\circ\Pcal)^*\circ\m_r^*(\Ical_\Ocal)\\
  &\subseteq(\m_m\circ\Pcal)^*(\Ical_\Ocal)
    &&\text{since $\m_r^*(\Ical_\Ocal)\subseteq\Ical_\Ocal$,}\\
  &=(m\Pcal)^*(\Ical_\Ocal) = D_{m\Pcal}\cdot\ZZ
    &&\text{by definition of $D_{m\Pcal}$.}
\end{alignat*}
The one point that possibly requires further explanation is the
inclusion $\m_r^*(\Ical_\Ocal)\subseteq\Ical_\Ocal$ of ideal sheaves
on~$\Gcal$.  The validity of this inclusion follows from the following
two facts:
\begin{itemize}
\item
The sheaf~$\Ical_\Ocal$ is the ideal sheaf of the image~$\Ocal(\ZZ)$ of
the identity section with its induced-reduced subscheme structure.
\item
The zero section satisfies $r\Ocal=\Ocal$, so
$\m_r(\Ocal(\ZZ))=(r\Ocal)(\ZZ)=\Ocal(\ZZ)$ as subsets of~$\Gcal$.
\end{itemize}
This proves that $D_{n\Pcal}\cdot\ZZ\subseteq D_{m\Pcal}\cdot\ZZ$, which is
equivalent to~$D_{m\Pcal}|D_{n\Pcal}$.
\end{proof}

\begin{definition}
A \emph{geometric divisibility sequence} is the divisibility
sequence~$(D_{n\Pcal})_{n\ge1}$ associated to a point
(section)~$\Pcal$ of infinite order in a group scheme~$\Gcal/\ZZ$ as
in Proposition~\ref{proposition:Dnisdivseq}
\end{definition}

In some cases an algebraic group~$G/\QQ$ has a particularly nice model
over~$\ZZ$, as for example is the case for abelian varieties. This
prompts the following definition.

\begin{definition}
Let $A/\QQ$ be an abelian variety and let $P\in A(\QQ)$ be a point of
infinite order. The \emph{abelian divisibility sequence} associated
to~$P$ is the divisibility sequence associated to the lift~$\Pcal$
of~$P$ to a section of the N\'eron model~$\Acal/\ZZ$ of~$A/\QQ$.
By abuse of notation, we denote this sequence by~$(D_{nP})_{n\ge1}$.
\end{definition}

We next show that Vojta's conjecture implies a strong upper bound for
abelian divisibility sequences on abelian vareities of dimension at
least~$2$. This result generalizes Theorem~\ref{theorem:ECoverQ} (take
\text{$A=E\times E$}).

\begin{proposition}
\label{proposition:hgcdonAV}
Let $A/\QQ$ be an abelian variety of dimension at least~$2$, and
assume that Vojta's conjecture is true for~$A$ blown up at~$O$. Fix a
Weil height
\begin{equation}
  \label{equation:htdefonA}
  h:A(\QQ)\to\RR
\end{equation}
on~$A$ with respect to an ample symmetric divisor.
\begin{parts}
\Part{(a)}
For every $\e>0$ there is a constant~$C=C(A,\e)$ and a proper
algebraic subvariety~\text{$Z\subsetneq A$} so that 
\[
  h_{\gcd}(P;O) \le \e h(P) + C
  \qquad\text{for all $P\in A(\QQ)\setminus Z$.}
\]
The exceptional set~$Z$ consists of a finite union of translates of
nontrivial abelian subvarieties of~$A$, so in particular, if~$A$ is
simple, then we may take~$Z=\emptyset$.
\Part{(b)}
Let~$(D_{nP})_{n\ge1}$ be the abelian divisibility
sequence associated to a point of infinte order~$P\in A(\QQ)$, and
assume further that the group~$\ZZ P$ generated by~$P$ is Zariski
dense in~$A$. Then for every~$\e>0$ there is a constant $C=C(A,P,\e)$
so that
\[
  \log D_{nP} \le \e n^2 + C
  \qquad\text{for all $n\ge1$.}
\]
\end{parts}
\end{proposition}

\begin{remark}
We observe that Proposition~\ref{proposition:hgcdonAV} is false if~$A$
is an elliptic curve, since then we have
$h_{\gcd}(P;O)=h_{E,O}(P)$ and $\log D_{nP}\sim n^2\hhat(P)$.  The
reason that our proof of Proposition~\ref{proposition:hgcdonAV} fails
when $\dim(A)=1$ is the requirement in
Theorem~\ref{theorem:vojta=>hgcdbound} that the subvariety~$Y$ have
codimension at least~2 in~$X$.
\end{remark}

\begin{proof}[Proof of Proposition~\ref{proposition:hgcdonAV}]
(a)
We apply Theorem~\ref{theorem:vojta=>hgcdbound} to the variety~$A$,
the subvariety consisting of the single point~$O$, and the ample
divisor used to define the height~\eqref{equation:htdefonA}.  The
canonical divisor on~$A$ is trivial, so
Theorem~\ref{theorem:vojta=>hgcdbound} says that
there is a subvariety $Z\subsetneq A$ such that
\[
  h_{\gcd}(P;O) \le \e h(P) + O(1)
  \qquad\text{for all $P\in A(\QQ)\setminus Z$.}
\]
This proves~(a), other than the characterization of~$Z$.
Let~$Z'\subset Z$ be any irreducible subvariety of~$Z$. If~$Z'(\QQ)$ is
finite, then we may discard it and adjust the~$O(1)$ accordingly. And
if~$Z'(\QQ)$ is infinite, then Faltings' theorem~\cite{Faltings} says that~$Z'$
is a translate of an abelian subvariety of~$A$.
\par(b)
We compute
\begin{alignat*}{2}
  \log D_{nP} &\le h_{\gcd}(nP;O) + O(1) 
    &\quad&\text{from Proposition~\ref{proposition:logDPlehgcd},} \\
  &\le \e h(nP) + O(1)
    &&\text{from (a), assuming $nP\notin Z$,} \\
  &\le \e n^2 \hhat(P) + O(1)
    &&\text{canonical height property \cite[B.5.1]{HS}}
\end{alignat*}
The fact that~$P$ has infinite order implies that~$\hhat(P)>0$, so
after replacing~$\e$ with~$\e/\hhat(P)$, this completes the proof
of~(b) provided \text{$nP\notin Z$}.
\par
Suppose that~$Z\ne\emptyset$, and let~$Z_1$
be an  irreducible component of~$Z$ that contains infinitely many
multiples of~$P$. From~(a), we know that
$Z_1=A_1+R_1$ for an abelian subvariety $A_1\subsetneq A$ and a point
$R_1\in A(\QQ)$. Choose $n_2>n_1$ with $n_1P\in Z_1$ and $n_2P\in Z_1$.
Then $(n_2-n_1)P\in A_1$. Letting $N=n_2-n_1$, it follows
that $P\in A_1+A[N]$, and hence that $nP\in A_1+A[N]$ for all $n\ge1$.
This contradicts the assumption that~$\ZZ P$ is Zariski dense in~$A$,
and hence there is no exceptional set.
\end{proof}

\section{Final remarks and questions}
\label{section:finalremarks}

We have proven a number of strong bounds for generalized greatest
common divisors and divisibility sequences, all conditional on 
the validity of Vojta's beautiful, but deep, conjecture applied to
an appropriate blowup variety. It would be of great interest to
find unconditional proofs of some of these results.
\par
In addition to height bounds, there are many other natural questions
that one might ask about abelian, or more generally geometric,
divisibility sequences.  For example, which such sequences contain
infinitely many prime numbers (cf.~\cite{EW1}). This is, of course, a
notoriously difficult question, even for the simplest divisibility
sequence \text{$2^n-1$}. There is some evidence~\cite{EGW} that
elliptic divisibility sequences~$(D_{nP})_{n\ge1}$ donot contain
infinitely many primes, although more general elliptic divisibility
``sequences'' $(D_{nP+mQ})_{m,n\ge1}$ may well contain infinitely many
primes.
\par
One might ask if a geometric divisibility sequence necessarily 
grows, or if it often returns to small values. 
For example, Ailon and Rudnick~\cite{AR} conjecture that if~$a,b\in\ZZ$
are multiplicatively independent, then
\[
  \gcd(a^n-1,b^n-1)=\gcd(a-1,b-1)
  \qquad\text{for infinitely many $n\ge1$.}
\]
They prove a strong version of this with~$\ZZ$ replaced by the polynomial
ring~$\CC[T]$. (See also~\cite{JHSFF} and~\cite{JHSECFF}
for analogs over~$\FF_q[T]$ and for elliptic curves.) 
We certainly suspect that the same is true for semiabelian varieties.

\begin{conjecture}
\label{conjecture:DnP=DPinftyoften}
Let~$\Gcal/\ZZ$ be a group scheme, let $\Pcal\in\Gcal(\ZZ)$ be
a~$\ZZ$-valued point, and assume that the following are true:
\begin{parts}
\Part{(1)}
The generic fiber~$G=\Gcal\times_\ZZ\QQ$ is an irreducible 
commutative algebraic group of dimension at least~$2$ with no unipotent part.
\Part{(2)}
Let~$P\in G(\QQ)$  be the restriction of~$\Pcal$ to the generic fiber.
Then the subgroup $\ZZ P$ generated by~$P$ is Zariski dense in~$G$.
\end{parts}
Then the geometric divisibility sequence $(D_{n\Pcal})_{n\ge1}$
corresponding to~$\Pcal$ satisfies
\[
  D_{n\Pcal} = D_{\Pcal}
  \qquad\text{for infinitely many $n\ge1$.}
\]
\end{conjecture}

It is tempting to guess that something similar is true for geometric
divisibility sequences associated to any irreducible algebraic group
of dimension at least~$2$, regardless of whether or not it is
commutative. (Note that the Zariski density condition is vital.) But
with no significant evidence for even
Conjecture~\ref{conjecture:DnP=DPinftyoften}, we will be content to
leave the general case as a question.


\end{document}